\documentclass[12pt]{amsart}
\usepackage{amsmath,amsthm,amsfonts,amssymb}
\begin{document} 
\newcommand{\Alb}{\mathop{\mathrm{Alb}}}
\newcommand{\m}{\mathfrak{m}}
\newcommand{\Spec}{{\mathrm{Spec}~}}
\newcommand{\Sing}{\mathop{\mathrm {Sing}}}
\newcommand{\supp}{\mathop{supp}}
\newcommand{\tensor}{\otimes}
\newcommand{\codim}{\mathop{codim}}
\newcommand{\B}{{\mathbb B}}
\newcommand{\C}{{\mathbb C}}
\newcommand{\N}{{\mathbb N}}
\renewcommand{\O}{{\mathcal O}}
\newcommand{\Q}{{\mathbb Q}}
\newcommand{\Z}{{\mathbb Z}}
\renewcommand{\P}{{\mathbb P}}
\newcommand{\R}{{\mathbb R}}
\newcommand{\rc}{\subset}
\newcommand{\rank}{\mathop{rank}}
\newcommand{\trace}{\mathop{tr}}
\newcommand{\dimc}{\mathop{dim}_{\C}}
\newcommand{\Lie}{\mathop{Lie}}
\newcommand{\Aut}{\mathop{{\rm Aut}}}
\newcommand{\Auto}{\mathop{{\rm Aut}_{\mathcal O}}}
\newcommand{\alg}[1]{{\mathbf #1}}
\newtheorem{corollary}{Corollary}
\newtheorem*{lemma}{Lemma}
\newtheorem*{proposition}{Proposition}
\newtheorem{theorem}{Theorem}
\newtheorem*{MainTheorem}{Main Theorem}
\theoremstyle{definition}
\newtheorem*{Beweis}{{\usefont{\encodingdefault}{\familydefault}{\seriesdefault}%
{\shapedefault} Proof}}
\newtheorem{remark}{{\usefont{\encodingdefault}{\familydefault}{\seriesdefault}%
{\shapedefault}\scshape Remark}}
\title{%
On Manifolds with trivial logarithmic tangent bundle
}
\author {J\"org Winkelmann}
%
\address{%
J\"org Winkelmann \\
 Institut Elie Cartan (Math\'ematiques)\\
 Universit\'e Henri Poincar\'e Nancy 1\\
 B.P. 239, \\
 F-54506 Vand\oe uvre-les-Nancy Cedex,\\
 France
}
\email{jwinkel@member.ams.org\newline\indent{\itshape Webpage: }%
http://www.math.unibas.ch/\~{ }winkel/
}
\thanks{
{\em Acknowledgement.}
The author wants to thank 
the Korea Institute for Advanced Study (KIAS) in Seoul.
The research for this article was done during the stay of the
author at this institute.}
\maketitle
\section{Introduction}
By a classical result of Wang \cite{W} a connected
compact complex manifold $X$
has holomorphically trivial tangent bundle if and only if there is
a connected complex Lie group $G$ and a discrete subgroup $\Gamma$ such that
$X$ is biholomorphic to the quotient manifold $G/\Gamma$.
In particular $X$ is homogeneous.
If $X$ is K\"ahler, $G$ must be commutative and the quotient manifold
$G/\Gamma$ is a compact complex torus.

The purpose of this note is to generalize this result to the
non-compact K\"ahler case. 
Evidently, for arbitrary non-compact complex manifold
such a result can not hold. For instance, every domain over $\C^n$
has trivial tangent bundle, but many domains have no automorphisms.

So we consider the {\em ``open case''} in the sense of Iitaka (\cite{I}),
i.e. we consider manifolds which can be compactified
by adding a divisor.

Following a suggestion of the referee, instead of only considering K\"ahler
manifolds we consider manifolds in class ${\mathcal C}$
as introduced in \cite{Fj}.
A compact complex manifold $X$ is said to be class in ${\mathcal C}$
if there is a surjective holomorphic map from a compact K\"ahler manifold
onto $X$. Equivalently, $X$ is bimeromorphic to a K\"ahler manifold
(\cite{V}).
For example, every Moishezon manifold is in class ${\mathcal C}$.

We obtain the following characterization:
\begin{MainTheorem}
Let $\bar X$ be a connected compact complex manifold, $D$ a closed analytic
subset and $X=\bar X\setminus D$.
Assume that $X$ is in class ${\mathcal C}$ 
as introduced in \cite{Fj} 
(also called ``weakly K\"ahler'').

Then the following conditions are equivalent:
\begin{enumerate}
\item
$D$ is a divisor with only simple normal crossings as singularities
and
the logarithmic tangent bundle $T(-\log D)$ is 
a holomorphically trivial vector bundle on $\bar X$.
\item
There is a complex semi-torus $G$ acting effectively on $\bar X$ with
$X$ as open orbit such that the all the isotropy groups are 
themselves semi-tori.
\end{enumerate}

Moreover, if one (hence both) of these conditions
are fulfilled, then there is a short exact sequence of
complex Lie groups
\[
1 \to (\C^*)^d\to G \to \Alb(\bar X) \to 1
\]
where $\Alb(\bar X)$ denotes the Albanese torus of $\bar X$
and $d=\dim(X)-\frac{1}{2}b_1(X)$.
\end{MainTheorem}

In the algebraic category we have the following result.

\begin{corollary}
Let $\bar X$ be a non-singular complete algebraic variety defined over $\C$, 
$D$ a divisor with only
simple normal crossings as singularities and let $X=\bar X\setminus D$.

Then the following conditions are equivalent:
\begin{enumerate}
\item
The logarithmic tangent bundle $T(-\log D)$ is 
a trivial vector bundle on $\bar X$.
\item
There is a semi-abelian variety $T$ acting on $\bar X$ with
$X$ as open orbit.
\end{enumerate}
\end{corollary}

\begin{corollary}
Let $X$ be a nonsingular algebraic variety defined over $\C$.

Then the following are equivalent:
\begin{enumerate}
\item
There exists a completion $X\hookrightarrow\bar X$ such that 
$D=\bar X\setminus X$ is a s.n.c.{} divisor and $T(-\log D)$ is trivial.
\item
$X$ is isomorphic to a semi-abelian variety.
\end{enumerate}
\end{corollary}
\begin{Beweis}
This follows from the preceding result because every semi-abelian variety
admits an equivariant completion by a s.n.c.{} divisor.
\end{Beweis}

\section{Terminology}
A {\em complex semi-torus} $G$ is a connected complex Lie group which can
be realized as a quotient of a vector group $(\C^n,+)$ by a discrete
subgroup $\Gamma$ such that $\Gamma$ generates $\C^n$ as
complex vector space.
Equivalently, a connected complex Lie group $G$ is a {\em semi-torus}
if and only if there exists
a short exist sequence of complex Lie groups
\[
0 \to L \to G \to T \to 0
\]
where $T$ is a compact complex torus and $L\simeq(\C^*)^k$ for some
$k\in\N$.

A divisor $D$ on a complex manifold $X$ has only 
{\em simple normal crossings}
as singularities if for every point $x\in X$ there exists local
coordinates $z_1,\ldots,z_n$ and a number $k\in\{0,\ldots n\}$ such that
in a neighbourhood of $x$ the divisor $D$ equals the zero divisor
of the holomorphic function $\Pi_{i=1}^k z_i$.
For a brevity, such a divisor is called a {``\em s.n.c.''} divisor.


Let $\bar X$ be a compact complex manifold with s.n.c.{} divisor $D$.
There is a stratification as follows: $X_0=X=\bar X\setminus D$,
$X_1=D\setminus \Sing(D)$ and for $k>1$ the stratum $X_k$ is the
non-singular part of $\Sing(\bar X_{k-1})$.
If in local coordinates $D$ can be written as $\{z:\Pi_{i=1}^d z_i=0\}$,
then $z=(z_1,\ldots,z_n)\in X_k$ iff $\#\{i:1\le i \le d, z_i=0\}=k$.

Let $D$ be an effective divisor on a complex manifold $\bar X$.
Then the sheaf $\Omega^1(\log D)$ of 
{\em logarithmic $1$-forms} with respect to $D$ is defined as
the $\O_X$-module subsheaf of the sheaf of meromorphic one-forms on
$\bar X$ which is locally generated by all
$df/f$ where $f$ is a section ${\mathcal O}_{\bar X}\cap{\mathcal O}_X^*$.

This sheaf is always coherent. It is locally free if
$D$ is a s.n.c.{} divisor. In fact, if
$D=\{z_1\cdot\ldots\cdot z_d=0\}$, then $\Omega^1(\log D)$
is locally the free $\O_X$-module over 
$dz_1/z_1,\ldots,dz_d/z_d,dz_{d+1},\ldots,dz_n$.

For a {\em s.n.c.} divisor $D$ on a complex manifold
we define the {\em logarithmic
tangent bundle} $T(-\log D)$ as the dual bundle of $\Omega^1(\log D)$.

Then $T(-\log D)$ can be identified
with the sheaf of those holomorphic vector fields $V$ on $\bar X$ 
which fulfill the following property:

{\em $V_x$ is tangent to $X_k$ at $x$ for every $k$ and 
every $x\in X_k$}.

In local coordinates: If $D=\{z:\Pi_{i=1}^d z_i=0\}$, then
$T(-\log D)$ is the locally free sheaf generated by
the vector fields $z_i\frac{\partial}{\partial z_i}$ ($1\le i\le d$)
and $\frac{\partial}{\partial z_i}$ ($d<i\le n$).

\section{The Proof of the Main Theorem}

\begin{Beweis}
$(1)\Rightarrow(2):$

Triviality of $T(-\log D)$ implies that the sheaf of logarithmic
one-forms $\Omega^1(\log D)$ is trivial as well.

Let $V=\Omega^1(\bar X,\log D)$ and $V^*=\Gamma(\bar X,T(-\log D))$. 
By \cite{D},\cite{N} every logarithmic one-form 
$\omega\in\Omega^1(\bar X,\log D)$
is closed if $X$ is K\"ahler. For an arbitrary manifold $X$ in class
${\mathcal C}$ there is always a holomorphic surjective bimeromorphic
map
$p:X'\to \bar X$ from some compact K\"ahler manifold $X'$. Moreover,
by blowing-up $X'$ if necessary, we may assume that $p^{-1}(D)$
is a s.n.c~divisor on $X'$. Now the aforementioned result for K\"ahler
manifolds implies that $d(p^*\omega)=0$ for every $\omega
\in\Omega^1(\bar X,\log D)$.
Since $p$ is biholomorphic on some open subset, we obtain $d\omega=0$.
Therefore closedness of logarithmic one-forms holds not only for K\"ahler
compact complex manifolds, it holds for all manifolds in class 
${\mathcal C}$.

Thus
\[
0=d\omega(x,y)= x(\omega(y))-y(\omega(x))-\omega([x,y])
\]
for $\omega\in V$, $x,y\in V^*$. 
Now $\omega(y)$ and $\omega(x)$ are global holomorphic functions
on a compact manifold and therefore constant.
Hence $x(\omega(y))=0=y(\omega(x))$. It follows that
$\omega([x,y])=0$ for all $\omega\in V$, $x,y\in V^*$. 
Thus $V^*$ is a {\em commutative} Lie algebra
of holomorphic vector fields on $\bar X$.
Let $G\subset\Aut(\bar X)$ denote the subgroup generated by the one-parameter
groups corresponding to vector fields $v\in V^*$. 
Recall that the sections in $T(-\log D)$ can be regarded as the vector fields
which are tangent to $X_k$ at every $x\in X_k$ for all $k$.
It follows that the $G$-orbits in $\bar X$ are precisely the connected
components of the strata $X_k$.
In particular $G$ has an open orbit, namely $X=X_0=\bar X\setminus D$.
Furthermore the closed orbits of $G$ are the connected components
of the unique closed stratum $X_d$ 
where $d$ is the largest natural number with $X_d\ne\emptyset$.

The existence of an open orbit implies that $G$ acts transitively 
on the Albanese torus $\Alb(\bar X)$.
Therefore, all the fibers of $\bar X\to\Alb(\bar X)$ are isomorphic.
Let $\bar X\to Y\to \Alb(\bar X)$ be the Stein factorization.
Since the Stein factorization is canonical, it is compatible with the
$\Aut(\bar X)$-action. 
For this reason $Y\to\Alb(X)$ is a finite holomorphic
map of $G$-spaces. Hence $G$ acts transitively on $Y$. It follows
that $Y$ is a compact complex space which is a quotient of a connected
commutative complex Lie group, in other words, $Y$ must be a compact
complex torus. By the universality property of the Albanese torus
this implies $Y=\Alb(\bar X)$.

Thus the fibers of $\bar\rho:\bar X\to\Alb(\bar X)$ are connected.
Let $H$ be the kernel of $G\mapsto\Alb(\bar X)$.
Recall that $G$ is commutative.
It follows that the $H$-orbits in $\bar X$
are precisely the intersections of $G$-orbits in $\bar X$
with fibers of the map $\bar\rho:\bar X\to\Alb(\bar X)$.
Moreover, $H$ acts freely on an open orbit in each fiber of
$\bar\rho$. This implies that $H$ is connected.

Now let $Z$ be a closed $G$-orbit (i.e. a connected component
of the smallest stratum $X_d$). Then the fibers of $\bar\rho|_Z$ 
are closed $H$-orbits.

If $\bar X$ is K\"ahler, a result of Sommese (\cite{S}, prop. 1)
implies that closed orbits of $H$ are fixed points.
Due to Fujiki (\cite{Fj}) the same assertion holds for
an arbitrary manifold $\bar X$ in class $\mathcal C$.

It follows that $\bar\rho|_Z:Z\to\Alb(\bar X)$ is biholomorphic.

For each irreducible component $D_i$ of $D$ and each $\omega\in V$
we may regard the {\em residue} $res_i(\omega)$. This residue is given
as integral of $\omega$ over a small loop around $D_i$.
A priori, it is a holomorphic function on $D_i$. But, since $D_i$
is compact, this holomorphic function is constant
(Alternatively, one may also use a calculation in local coordinates
which shows that $d\omega=0$ forces $res_i(\omega)$ to be locally
constant).

Let $n=\dim\bar X$ and $g=\dim\Alb(\bar X)$.

Fix $p\in Z$.
Near $p$, $Z$ is the intersection of $d=n-g$ irreducible components
$D_1,\ldots,D_d$ of $D$. We choose a basis $(\omega_1,\ldots,\omega_n)$
of $V$ such that
\[
res_i(\omega_j)_p=
\begin{cases}
2\pi i & \text{ if $i=j\le d$},\\
0 & \text{ if  $i\ne j$ or $j>d$}
\end{cases}
\]
Now we choose a point $q\in X$ near $p$ and
define local coordinates $z_i$ near $p$ via
\[
z_i(x)=
\begin{cases}
\exp\left(\int_q^x \omega_i\right) &  \text{ if $i\le d$},\\
\int_p^x \omega_i & \text{ if  $i>d$}
\end{cases}
\]
Then in these local coordinates we can describe the $\omega_i$ as follows:
\[
\omega_i=
\begin{cases}
\frac{dz_i}{z_i} & \text{if $i\le d$}\\
dz_i & \text{ if $i>d$ }
\end{cases}
\]
It follows that there is a biholomorphic map from a neighbourhood of
$p$ to a neighbourhood of $0$ in $\C^n$ taking the fundamental vector fields
of $\Lie(H)$ into vector fields of the form $\sum_{i\le d} a_iw_i\frac{\partial}%
{\partial w_i}$.
This implies in particular that $H$ contains a totally real compact subgroup 
$K=(S^1)^d$ acting as
\[
K\ni(\theta_1,\ldots,\theta_d):z\mapsto(\theta_1z_1,\ldots,\theta_dz_d,
z_{d+1},\ldots,z_n)
\]
Thus $H$ is a connected commutative
complex Lie group of dimension $d$ containing a 
totally real compact subgroup of real dimension $d$. It follows that
$H$ is a complex semi-torus.
On the other hand, the above description of the $H$-vector fields 
in local coordinates also
implies that $H$ admits an almost faithful representation
on the tangent space of each fixed point of $H$.
Therefore $H$ is a semi-torus of dimension $d$ which admits an almost faithful
representation. It follows that 
$H$ must be isomorphic to $(\C^*)^d$.

As a consequence we obtain that $G$ admits a short exact sequence
of complex Lie groups in the form
\[
1 \to (\C^*)^d \to G \to \Alb(\bar X) \to 1
\]
Thus $G$ is a semi-torus as well.

$(1)\Leftarrow(2):$

Let $v_1,\ldots,v_n$ be a basis for the vector space of $G$-fundamental
vector fields on $\bar X$. The complement $D=\bar X\setminus X$ of the
open orbit $X$ can be characterized as the set of those points where
the vector fields $v_i$ fail to span the tangent bundle $T_{\bar X}$.
Thus $D$ is defined by the vanishing of 
\[
w=\Lambda_{i=1}^n v_i\in\Gamma\left(\bar X,\Lambda^nT_{\bar X}\right)
\]
Since $\Lambda^nT_{\bar X}$ is a line bundle, it is clear that $D$ is
of pure codimension one.

Now let $p\in \supp D$ and $G_p=\{a\in G:a(p)=p\}$.
By assumption, $G_p$ is a semi-torus and therefore reductive. 
This implies that the $G_p$-action
can be linearized near $p$,  i.e. there is an $G_p$-equivariant
biholomorphism between an open neighbourhood of $p$
in $\bar X$ and an open neighbourhood of $0$ in the vector space
$W=T_p\bar X$.
It follows that a neighbourhood of $p$ in $D$ is isomorphic
to a union of vector subspaces of codimension one in $W$. From
the assumption that the $G_p$-action is effective, one can
deduce that this is a transversal union, i.e.~$D$ is a simple normal
crossings divisor near $p$.

Recall that the isotropy groups are required to be semi-tori.
In particular, they are reductive.
Therefore
the action of every isotropy subgroup is linearizable in some
open neighbourhood.
This implies that for every point on $\bar X$ we can find
a system of local coordinates in which the $G$-fundamental
vector fields are simply
\[
\frac{\partial}{\partial z_1},\ldots,\frac{\partial}{\partial z_d},
z_{d+1}\frac{\partial}{\partial z_{k+1}},\ldots,z_{n}\frac{\partial}{\partial z
_n},
\]
where $d$ equals the co
dimension of the isotropy group.
Hence we have
\begin{equation*}
{\bf T}(\bar X, D)\cong \bar X \times \Lie G.
\end{equation*}
\end{Beweis}

\begin{remark}\label{thm-nonkaehler}
In the proof for the direction ``$(1)\Leftarrow(2)$'' we did not employ
the K\"ahler assumption. Therefore this part of the theorem
is valid even without requiring $X$ to be K\"ahler.
\end{remark}

\section{Examples}
\subsection{Toric varieties}
The easiest examples of equivariant compactifications of $(\C^*)^d$
with trivial logarithmic tangent bundles are $\P_d(\C)$ and $\P_1(\C)^d$.
More examples are obtained from the theory of toric varieties,
see e.g. \cite{Fu},\cite{O}.

Now let $A$ be a complex semi-torus admitting a short exact sequence
\begin{equation}\label{ses}
1 \to (\C^*)^d \to A \to T \to 1
\end{equation}
where $T$ is a compact complex torus.
Let $L\hookrightarrow\bar L$ be a smooth equivariant
compactification of $L=(\C^*)^d$. Then a smooth equivariant
compactification of $A$ can be constructed as a fiber product:
$\bar A=(A\times\bar L)/\sim$ where $(a,x)\sim(a',x')$
iff there exists an element $g\in L$ such that $a\cdot g^{-1}=a'$
and $g\cdot x=x'$.

This construction preserves the K\"ahler condition:
\begin{lemma}
If $\bar L$ is K\"ahler, then $\bar A$ is K\"ahler, too.
\end{lemma}
\begin{Beweis}
The fiber bundle $(\ref{ses})$ is given by locally constant
transition functions with values in the maximal compact subgroup
$K$ of $L$, acting by multiplication. By averaging we may assume
that the K\"ahler metric on $\bar L$ is $K$-equivariant.
Thus the associated K\"ahler form $\omega_1$ induces a closed
semi-positive $(1,1)$-form on $\bar A$ such that the restriction
to the tangent bundle of any fiber is positive.
Taking sum of this $(1,1)$-form and the pull-back of a K\"ahler
form on $T$ yields a K\"ahler form on $\bar A$.
\end{Beweis}

(This is a special case of a general result of Blanchard 
\cite{B1} which implies
that for any holomorphic fiber bundle of compact complex
manifolds $E\to B$ with typical fiber
$F$ and $b_1(F)=0$ the K\"ahler property for both $B$ and $F$ implies
the K\"ahler property for $E$.)

In this way we see that every semi-torus admits a smooth
equivariant K\"ahler compactification. On the other hand, our main
theorem implies that every smooth equivariant K\"ahler compactification
of a semi-torus arises in this way.

In contrast, non-K\"ahler compactifications may arise in many ways,
see e.g.~\cite{LM} and the examples given further below in this
article.
\subsection{Non-K\"ahler examples}
If $X$ is a compact complex manifold with trivial tangent bundle,
and $G$ the connected component of its automorphism group, then
$X$ is K\"ahler if and only if $G$ is commutative.

In the logarithmic case, there is such a conclusion only in one direction:
If $X$ is a compact complex manifold with s.n.c.{} divisor $D$ such that
$T(-\log D)$ is trivial, then the K\"ahler assumption implies that the
connected component of $\Aut(X,D)$ is commutative.
On the other hand the commutativity does not imply the K\"ahler property
as we will see by the example given below.

Let $\alpha,\beta\in\C$ with $|\alpha|,|\beta|>1$.
We define a $\Z$-action on $\C^2\setminus\{(0,0)\}$ by
\[
(z_1,z_2)\mapsto (\alpha^nz_1,\beta^nz_2)
\]
Then the quotient of $\C^2\setminus\{(0,0)\}$ by this $\Z$-action
is a so-called {\em Hopf surface}. Such a Hopf surface $\bar X$ is 
diffeomorphic
to $S^1\times S^3$. In particular $\dim H^1(\bar X,\C)$ is odd and 
therefore $\bar X$ can not be K\"ahler.

Now let $T$ be the quotient of $\C^*\times\C^*$ by the subgroup
\[
\{(\alpha^n,\beta^n):n\in\Z\}.
\]
Then $T$ is a complex semi-torus and $\bar X$ is an equivariant
compactification of $T$. 
The isotropy groups at the two non-open orbits
are isomorphic to $\C^*$. 

Thus all the isotropy groups are semi-tori and consequently
the logarithmic tangent bundle
is trivial (see Remark \ref{thm-nonkaehler}).

\subsection{The noncommutative case}
Let $\bar X$ be a non-K\"ahler compact complex manifold
with s.n.c.{} divisor $D$ such that $T(-\log D)$ is trivial.

In this case it is still true that there is a connected complex Lie group
$G$ with $\dim(G)=\dim(\bar X)$ acting on $\bar X$ with $X=\bar X\setminus D$
as open orbit. However, $G$ might be non-commutative and the $G$-action on
the open orbit is only almost free, i.e., the isotropy group at a point
of the open orbit is not necessarily trivial, but at least discrete.

The easiest such examples, with $D=\emptyset$, are obtained by
considering discrete subgroups $\Gamma$ in connected complex Lie groups
with compact quotient $G/\Gamma$. By a result of Borel 
(see \cite{B}) every semisimple
Lie group contains such a ``cocompact'' discrete subgroup $\Gamma$.
Such complex quotients have been studied in \cite{Wi}.

Next let us give an example with $D\ne\emptyset$.
Recall that
$SL_2(\C)$ contains discrete cocompact subgroups $\Gamma$ such that
$b_1(Y)>0$ for $Y=SL_2(\C)/\Gamma$ (see \cite{J}). Then $H^1(Y,{\mathcal O})$
is a vector space of positive dimension equal to $b_1(Y)>0$ and
the induced action of $SL_2(\C)$ on $H^1(Y,{\mathcal O})$ is trivial
(see \cite{Akh}).
For any $\alpha\in H^1(Y,{\mathcal O})$ let $\alpha'$ denote the
image via 
\[
\exp:H^1(Y,{\mathcal O})\to H^1(Y,{\mathcal O}^*).
\]
Then $\alpha'$ defines a topologically trivial $\C^*$-principal
bundle. Since the $SL_2(\C)$-action on $H^1(Y,{\mathcal O})$ is trivial,
this is a homogeneous bundle. If $G$ denotes the connected component
of the group of $\C^*$-principal bundle automorphisms, we thus obtain
a short exact sequence
\[
1 \to \C^* \to G \to SL_2(\C) \to 1.
\]
Such a sequence is necessary split. Hence $G\simeq\C^*\times SL_2(\C)$.
Now consider the compactification $\bar X$ of the total space $X$ of this
bundle given by adding a $0$- and a $\infty-section$.
Then $\bar X$ has a trivial logarithmic tangent bundle $T(-\log D)$
for $D=\bar X\setminus X$ and $G$ acts on $\bar X$ with $X$ as open orbit.

\subsection{The condition on the isotropy groups}
Next we present an example to show
that the condition on the isotropy groups in property $(2)$ in the
Main theorem can not be dropped.

To see this, we first note that a complex-analytic semi-torus 
may contain closed complex Lie subgroups which are not semi-tori.
The easiest such example is obtained as follows:
Embedd the additive group (which is not a semi-torus)
into the semi-torus $\C^*\times\C^*$ via
\[
\C\ni t\stackrel{f}{\mapsto} (e^t,e^{it})\in\C^*\times\C^*
\]
To verify that the image is indeed a {\em closed} subgroup,
consider its pre-image in the universal covering of $\C^*\times\C^*$.
If we realize the universal covering by 
$\pi:(z_1,z_2)\mapsto (e^{z_1},e^{z_2})$,
then 
\[
\pi^{-1}(f(\C))= \{(t,it):t\in\C\}+(2\pi i\Z)^2=\{(z_1,z_2):
\frac{z_2-iz_1}{2\pi i}\in \Z[i] \}.
\]
Hence $f(\C)$ is closed and $(\C^*\times\C^*)/f(\C)\simeq (\C/\Z[i])$.
In this way the complex manifold $X=\C^*\times\C^*$ can be realized
as a $\C$-principal bundle over the elliptic curve $E=\C/\Z[i]$.
The embedding $\C\hookrightarrow\P_1(\C)$ induces a compactification
$\bar X$ of $X$ by adding a $\infty$-section to the $\C$-principal
bundle $X\to E$. Thus we obtain an equivariant compactification
where the
isotropy group for a point at $\bar X\setminus X$ is isomorphic to
$(\C,+)$ and therefore not a semi-torus.
If we now look at the vector field corresponding to this
action of the additive group $(\C,+)$, we see that it vanishes
of order two at the $\infty$-section.
This implies that, regarded as a section in the logarithmic tangent
bundle, it does vanish at the $\infty$-section with multiplicity one.
In particular, the logarithmic tangent bundle admits a holomorphic
section which vanishes somewhere, but not everywhere.
Hence the logarithmic tangent bundle can not be trivial.
\begin{remark}
This example was first used by Serre for an entirely different
purpose: Via GAGA the $\C$-principal bundle over $E$ is algebraic.
In this way one obtains an exotic algebraic structure on the complex
manifold $\C^*\times\C^*$. This yields an example of an algebraic variety
which is not affine (in fact every regular function is constant), although
it is Stein as a complex manifold.
\end{remark}
\section{An Application}
Let $U\to \Delta$ be a family of projective manifolds. By our
characterization of log-parallelizable manifolds we obtain an
easy proof that the set of all $t\in \Delta$ for which $U_t$ is an
equivariant compactification of a semi-abelian variety forms a
{\em constructible} subset of $\Delta$.

\begin{proposition}
Let $\pi:U\to \Delta$ be a smooth proper connected surjective map
between K\"ahler complex manifolds and let $D$ be a hypersurface on
$U$ which does not contain any fiber of $\pi$.

Let $S$ denote the set of all $t\in \Delta$ for which the fiber 
$U_t=\pi^{-1}(\{t\})$ is a equivariant compactification of a 
semi-abelian variety with $U_t\setminus D$ as open orbit.
 
Then $S$ is {\em constructible}, i.e., there is a finite family
of pairs of closed analytic subsets $Y_i\subset Z_i\subset \Delta$ 
such that 
$S=\cup_{i}Z_i\setminus Y_i$.
\end{proposition}
\begin{Beweis}
By our theorem, $S$ coincides with the set of all $t\in\Delta$
with the property that $\pi^{-1}(p)\cap D$ is a s.n.c.~divisor and 
furthermore the logarithmic tangent bundle is trivial.

The fiber dimension $\dim_x(\pi^{-1}(\pi(x)))$ is Zariski
semicontinuous, because $\pi$ is proper.
Hence there is no loss in generality in assuming that
the fiber dimension is constant.
Let $r$ denote this fiber dimension.
Then for every $p\in\Delta$ the number of irreducible components
of $\pi^{-1}(p)\cap D$ equals $\dim H^{2r-2}(\pi^{-1}(p)\cap D,\C)$.
Using the resolution of the sheaf of local constant functions
$\underline{\C}$ by coherent sheaves via the holomorphic
de Rham complex
\[
0 \to \underline{\C} \to \O \to \Omega^1 \to \ldots 
\]
combined with the semicontinuity results for coherent sheaves
it follows that $\Delta$ decomposes as a finite union of
constructible sets along which the number of irreducible
components of $\pi^{-1}(p)\cap D$ is constant.

Now a divisor with $m$ irreducible components in a compact complex
manifold $F$ is a {\em s.n.c.} divisor iff all the irreducible
components are smooth and meet transversally. These are Zariski
open conditions (as long as the number of irreducible
components does not jump).

Put together, these arguments show that the set $\Sigma$ of all $t\in\Delta$
for which the fiber $\pi^{-1}(p)\cap D$ is a s.n.c.~divisor
constitutes a constructible subset of $\Delta$.

Thus there is no loss of generality in assuming that $\pi^{-1}(p)\cap D$
is always a s.n.c.~divisor.

Let $E$ denote the bundle of logarithmic vertical vector fields.
Then $E$ is a vector bundle of rank $r$ where $r=\dim(U)-\dim(\Delta)$.
By our main theorem $t\in S$ iff $E|_{U_t}$ is holomorphically trivial.
By the semi-continuity theorem 
\[
Z=\left\{t\in \Delta: \dim \Gamma(U_t,E|_{U_t})\ge r 
\right\}
\]
is a closed analytic subset of $\Delta$
and
\[
\Omega=\left\{t\in \Delta: \dim \Gamma(U_t,E|_{U_t})=r\right\}
\]
is Zariski open in $Z$.
Moreover, for $t\in \Omega$ every section of $E|_{U_t}$  extends to some
neighbourhood of $U_t$ in $\pi^{-1}(\Omega)$.
Thus the set $W$ of all points in $\pi^{-1}(\Omega)$ where the sections
of $E|_{U_t}$ fail to span $E|_{U_t}$ is a closed analytic subset.
Since $\pi$ is proper, $S=\Omega\setminus\pi(W)$ is Zariski open in $Z$.
\end{Beweis}


\begin{thebibliography}{Bla}

\bibitem{Akh}
D.N.~Akhiezer: \sl
Group actions on the Dolbeault cohomology of homogeneous manifolds,
\rm  Math.~Z. \bf 226  \rm (1997), 607--621.

\bibitem{B1} A.~Blanchard: \sl
Sur les vari\'et\'es 
analytiques complexes,
\rm  Ann. Sci. Ecole Norm. Sup. \rm (3) {\bf 73} (1956), 157--202.

\bibitem{B}
A.~Borel: \sl
Compact Clifford Klein forms of symmetric spaces, 
\rm  Topology \bf 2 \rm (1963), 111--122.

\bibitem{D}
P.~Deligne: \sl
Theorie de Hodge. II,
\rm  Publ. Math. IHES \bf 40 \rm (1971), 5--57.

\bibitem{Fj}
A.~Fujiki:
\sl On automorphism groups of compact K\"ahler manifolds.
\rm Inventiones math. {\bf 44} (1978), 225--258.

\bibitem{Fu}
W.~Fulton: 
Introduction to Toric Varieties.
\rm Annals of Math.~Studies. Princeton University Press. {\bf 131},
1993.

\bibitem{I} 
S.~Iitaka: 
Algebraic Geometry. Springer GTM {\bf 76}, 1981.

\bibitem{J} 
T.~J\o rgensen:
\sl Compact 3--manifolds of constant negative curvature fibering over the
circle,
\rm  Ann. Math. \bf 106 \rm (1977), 61--72.

\bibitem{LM}
F.~Lescure, L.~Meersseman: \sl
Compactifications \'equivariantes non k\"ahl\'eriennes d'un groupe
alg\'ebrique multiplicatif,
\rm  Ann.~Inst.~Fourier \bf 52\rm, no. 1, (2002), 255--273.

\bibitem{N}
J.~Noguchi: \sl
A short analytic proof of closedness of logarithmic forms,
\rm  Kodai Math. J. \bf 18 \rm (1995), 295--299.

\bibitem{O}
T.~Oda: 
Convex Bodies and Algebraic Geometry.
Springer Verlag 1988.

\bibitem{S}
A.~Sommese: \sl
Holomorphic vector fields on compact K\"ahler manifolds,
\rm  Math.~Ann. \bf 210 \rm (1974), 75--82.

\bibitem{V}
J.~Varouchas: \sl
K\"ahler spaces and proper open morphisms. 
\rm Math.~Ann. \bf 283 \rm (1989), no. 1, 13--52.

\bibitem{W}
H.C.~Wang: \sl
Complex Parallisable manifolds,
\rm  Proc.~A.M.S. \bf 5 \rm (1954), 771--776.

\bibitem{Wi}
J.~Winkelmann: 
Complex-analytic geometry of complex parallelizable manifolds,
\rm  Memoirs de la S.M.F. \bf 72/73 \rm (1998).

\end{thebibliography}
\end{document}